\documentclass{amsart}

\usepackage{amssymb,amsmath,amsthm,latexsym,booktabs,todonotes}

\theoremstyle{definition}
\newtheorem{definition}{Definition}

\newtheorem{algorithm}[definition]{Algorithm}

\theoremstyle{plain}

\newtheorem{proposition}[definition]{Proposition}

\begin{document}

\title[Canonical forms of $2 \times 2 \times 2$ and $2 \times 2 \times 2 \times 2$ symmetric tensors]
{Canonical Forms of $2 \times 2 \times 2$ and $2 \times 2 \times 2 \times 2$ Symmetric Tensors over Prime Fields}

\author[Stavrou]{Stavros G. Stavrou}
\address{Department of Mathematics and Statistics, University of Saskatchewan, Canada}
\email{sgs715@mail.usask.ca}



\begin{abstract}
We consider symmetric tensors of format $2 \times 2 \times 2$ and $2 \times 2 \times 2 \times 2$ over prime fields. Using computer algebra we compute the canonical forms of these tensors. For $2 \times 2 \times 2$ symmetric tensors, we consider the prime fields $\mathbb{F}_p$ for $p = 2, 3, 5, 7, 11, 13, 17$. For $2 \times 2 \times 2 \times 2$ symmetric tensors, we consider the prime fields $\mathbb{F}_p$ for $p = 2, 3, 5, 7$.  For each canonical form, we determine the size of its orbit and the symmetric rank of the symmetric tensors in its orbit.
\end{abstract}

\maketitle

\section{Introduction}

A tensor is a multidimensional array of numbers. Formally, an order-$k$ tensor is an element of the tensor product of $k$ vector spaces. A first order tensor is a vector and a second order tensor is a matrix. These objects were studied by Cayley in the late 1800's and remain an active area of research \cite{Cayley} \cite{tenBerge}. The notion of tensor rank and the deomposition of a tensor into a sum of outer products of vectors was first introduced by Hitchcock \cite{Hitchcock1}  \cite{Hitchcock2}. Models of this decomposition were formed later: Harshman named his model PARAFAC (short for {\it parallel factor analysis}) \cite{Harshman}; and Carrol and Chang named their model CANDECOMP (short for {\it canonical decomposition}) \cite{CP}. Now these models are collectively referred to as CPD (short for {\it canonical polyadic decomposition}) \cite{KoldaBader}. 
For a textbook on tensor decompositions and related algorithms, as well as applications of tensors in areas such as data mining, email surveillance, gene expression classification, and signal processing, see Cichocki et al. \cite{Cichocki}, Kroonenberg \cite{Kroonenberg}, and Smilde et al. \cite{Smilde}. Symmetric tensors have applications in electrical engineering, signal processing, computational complexity, telecommunications, and data analysis; see  \cite{Bernardi}  \cite{ComonApplications1} \cite{ComonApplications2} \cite{Landsberg1} \cite{Landsberg2} and references therein. Bernardi et al. discuss algorithms for computing the symmetric rank of $2 \times \dots \times 2$ tensors from a geometric point of view \cite{Bernardi}. Landsberg et al. consider the rank and border rank of symmetric tensors using geometric methods \cite{Landsberg2}.

In this paper we use computer algebra to study symmetric tensors of format $2 \times 2 \times 2$ and $2 \times 2 \times 2 \times 2$ over prime fields, where we compute the canonical forms with respect to the action of the direct product of general linear groups GL$_{2,2,2}(\mathbb{F}_p)$ and GL$_{2,2,2,2}(\mathbb{F}_p)$, respectively.

More generally, the canonical forms of $2 \times 2 \times 2$ tensors have been studied many times in the literature \cite{Oldenburger2} \cite{Schwartz}. In 1881, Le Paige was the first to consider this problem over $\mathbb{C}$ \cite{LePaige}. Over $\mathbb{R}$, the problem was treated by Oldenburger \cite{Oldenburger1}. More recently, de Silva and Lim summarized these results \cite{deSilvaLim}. In recent work, we considered the problem of determining canonical forms of $2 \times 2 \times 2$ and $2 \times 2 \times 2 \times 2$ tensors over $\mathbb{F}_2$ and $\mathbb{F}_3$ \cite{BremnerStavrou}. 

For the $2 \times 2 \times 2 \times 2$ case, the problem of determining canonical forms over $\mathbb{R}$ and $\mathbb{C}$ has not yet been determined. This problem is related to quantum entanglement, where the classification of entanglements of four qubits corresponds to the determination of canonical forms of these tensors. See Luque and Thibon \cite{Luque}, Verstraete et al. \cite{Verstraete}, Borsten et al. \cite{Borsten}, Briand et al. \cite{Briand}, Chterental and Djokovic \cite{Chterental}, Brylinski \cite{Brylinski}, Gour and Wallach \cite{Gour}, Wallach \cite{Wallach}, Lamata et. al \cite{Lamata}, and references therein.  

For $2 \times 2 \times 2$ and $2 \times 2 \times 2 \times 2$ symmetric tensors over $\mathbb{R}$ and $\mathbb{C}$, the problem was solved by Gurevich \cite{Gurevich}. More recently, Weinberg \cite{Weinberg} independently reported the results again. In this paper, we extend our results in \cite{BremnerStavrou} by considering the same problem, but restricting our attention to symmetric tensors.

\section{Preliminaries}

\begin{definition}
An order-$k$ \textbf{tensor} $X$ is an element of the \textbf{tensor product} of $k$ vector spaces $V_1 \otimes V_2 \otimes \dots \otimes V_k$, where the \textbf{order} refers to the number of dimensions or modes. 
\end{definition}
This is an abstract definition that does not depend on a choice of basis in each vector space. This is one of the motivations for studying the canonical forms of tensors, since they represent properties of tensors which do not depend on the choice of basis. Once we fix a basis in each vector space, we can associate to each order-$k$ tensor a $k$-dimensional array. Then an order-1 tensor is a vector, and an order-2 tensor is a matrix, and we say that the order-$k$ tensor of format $d_1 \times \dots \times d_k$ is an element of $\mathbb{F}^{d_1 \times \dots \times d_k}$. In particular, we consider order-$k$ tensors $X = [x_{i_1 \dots i_k}]$ of format $2 \times \dots \times 2$ (with $k$ factors, for $k = 3, 4$) with entries in the prime field $\mathbb{F}_p$, where $i_1, \dots, i_k \in \{1, 2\}.$

Let's discuss different ways to represent tensors. One method is to write the entries of a tensor into vector (flattened) form. 

\begin{definition}\cite{KoldaBader}
The {\bf vectorization} or {\bf flattening} of $X$, denoted Vec$(X)$ writes the columns of $X$ in a vector format, where the entries are in lexical order of the $k$-tuples of the subscripts.
\end{definition}
We can also represent order-$k$, $k>2$ tensors $X$ in terms of matrix slices. We consider $2 \times 2 \times 2$ tensors.

\begin{definition}
A {\bf slice} of $X \in \mathbb{F}^{2 \times 2 \times 2}$ is a $2 \times 2$ matrix obtained by fixing one of the three indices. There are two slices in each of the three modes:
\begin{alignat*}{3}
& \text{{\bf frontal slices} (fix $k$):} 
\qquad 
& \begin{bmatrix}x_{111}&x_{121}\\x_{211}&x_{221}\end{bmatrix}, 
\qquad
& \begin{bmatrix}x_{112}&x_{122}\\x_{212}&x_{222}\end{bmatrix};
\\
& \text{{\bf vertical slices} (fix $j$):} 
\qquad 
& \begin{bmatrix}x_{111}&x_{112}\\x_{211}&x_{212}\end{bmatrix}, 
\qquad
& \begin{bmatrix}x_{121}&x_{122}\\x_{221}&x_{222}\end{bmatrix}; 
\\
& \text{{\bf horizontal slices} (fix $i$):} 
\qquad 
& \begin{bmatrix}x_{111}&x_{112}\\x_{121}&x_{122}\end{bmatrix}, 
\qquad
& \begin{bmatrix}x_{211}&x_{212}\\x_{221}&x_{222}\end{bmatrix}.
\end{alignat*}
\end{definition}
We conventionally represent $X \in \mathbb{F}^{2 \times 2 \times 2}$ in terms of its frontal slices:
\[
X
=
\left[
\begin{array}{cc|cc}
x_{111} & x_{121} & x_{112} & x_{122} \\
x_{211} & x_{221} & x_{212} & x_{222}
\end{array}
\right].
\]
Although we make use of representing $2 \times 2 \times 2$ tensors in terms of slices and its vectorized form, we restrict ourselves solely to the vectorized form of $2 \times 2 \times 2 \times 2$ tensors, since it's the simplest way to visualize them.

\begin{definition} \cite{Comon1}
An order-$k$ tensor $[x_{i_1 \dots i_k}] \in \mathbb{F}^{n \times \dots \times n}$ is called {\bf symmetric} if 
\[
x_{i_{\pi(1)} \dots i_{\pi(k)}} = x_{i_1 \dots i_k}, \qquad i_1, \dots, i_k \in \{1, \dots, n\},
\]
for all permutations $\pi \in S_k$. 
\end{definition}

\begin{proposition} [Proposition 3.7 \cite{Comon1}]
Let $X = [x_{i_1 \dots i_k}] \in \mathbb{F}^{n \times \dots \times n}$ be an order-$k$ tensor. Then 
\[
\pi(X) = X
\]
for all permutations $\pi \in S_k$ if and only if 
\[
x_{i_{\pi(1)} \dots i_{\pi(k)}} = x_{i_1 \dots i_k}, \qquad i_1, \dots, i_k \in \{1, \dots, n\}
\]
for all permutations $\pi \in S_k$.
\end{proposition}

For example, an order-3 tensor $[x_{ijk}] \in \mathbb{F}^{n \times n \times n}$ is symmetric if 
\[
x_{ijk} = x_{ikj} = x_{jik} = x_{jki} = x_{kij} = x_{kji}
\]
for all $i,j,k \in \{1, \dots, n\}$. In particular, consider the case when $n = 2$:
\[
X =
\left[
\begin{array}{cc|cc}
x_{111} & x_{121} & x_{112} & x_{122} \\
x_{211} & x_{221} & x_{212} & x_{222}
\end{array}
\right].
\]
Then $x_{112} = x_{121} = x_{211} =: b$ and $x_{122} = x_{212} = x_{221} =: c$, and so we can label $X$ in the following way:
\[
X =
\left[
\begin{array}{cc|cc}
a & b & b & c \\
b & c & c & d
\end{array}
\right].
\]

We will denote the set of all order-$k$, $n$-dimensional symmetric tensors over the field $\mathbb{F}$ by $\mathcal{S}^k(\mathbb{F}^n) \subset \mathbb{F}^{n \times \dots \times n}$. The set of such tensors satisfies the property $\pi(X) = X$ for all $\pi \in S_k$ and $X \in \mathcal{S}^k(\mathbb{F}^n)$.

\begin{definition}\cite{Hitchcock1} 
A tensor $X \in \mathbb{F}^{d_1 \times \dots \times d_k}$ is {\bf decomposable} if it can be written as
\[
X = u^{(1)} \otimes u^{(2)} \otimes \cdots \otimes u^{(k)},
\]
with non-zero $u^{(i)} \in \mathbb{F}^{d_i}$ for $i = 1, \dots, k$. The $(i_1, i_2, \dots, i_k)$th entry of $X$ is 
\[
x_{ i_1 i_2 \cdots i_k } = u^{(1)}_{i_1} u^{(2)}_{i_2} \cdots u^{(k)}_{i_k}.
\]
\end{definition}
A decomposable tensor is also called a {\bf simple} or {\bf rank-1} tensor. If we let $u^{(1)} = \dots = u^{(k)}$ in the definition above then we have defined a {\bf symmetric simple} tensor. This allows us to define two different notions of tensor rank.

\begin{definition}\cite{Hitchcock1}
A tensor has {\bf outer product rank} $r$ if it can be written as a sum of $r$ (and no fewer) decomposable tensors, 
\[
X = \sum_{i = 1}^r = u_i^{(1)} \otimes \dots \otimes u_i^{(k)} = u_1^{(1)} \otimes \dots \otimes u_1^{(k)} + \dots + u_r^{(1)} \otimes \dots \otimes u_r^{(k)}
\]
where $u_i^{(1)} \in \mathbb{F}^{d_1}, \dots, u_i^{(k)} \in \mathbb{F}^{d_k}, i = 1, \dots, r$.
We write rank$(X)$ to denote the outer product rank of $X$.
\end{definition}

\begin{definition}
A tensor has {\bf symmetric outer product rank} $s$ if it can be written as a sum of $s$ (and no fewer) symmetric simple tensors,
\[
X = \sum_{i=1}^s u_i^{\otimes k}
\]
We write rank$_S(X)$ to denote the symmetric outer product rank of $X$. 
\end{definition}
The only rank-0 (symmetric or not) tensor is the zero tensor. We will drop the words outer product, and simply say rank and symmetric rank. 

\begin{definition}
The {\bf maximum rank} is defined to be  
\[
\text{max}\{ \text{rank}(X) \mid X \in \mathbb{F}^{d_1 \times \dots \times d_k} \}. 
\]
If we replace rank$(X)$ with rank$_S(X)$ we get the analogous definition of maximum symmetric rank.
\end{definition}
Given a symmetric tensor $X$ we can compute its decomposition in $\mathcal{S}^k(\mathbb{F}^n)$ or in $\mathbb{F}^{n \times \dots \times n}$. Since the decomposition in the former set is constrained, it follows that for all $X \in \mathcal{S}^k(\mathbb{F}^n)$ we have the following inequality relating the two notions of rank:
\[
\text{rank}(X) \leq \text{rank}_S(X).
\]

In order to compute the canonical forms of symmetric tensors we use a tensor-matrix multiplication. {\bf Multilinear matrix multiplication} is a tensor-matrix multiplication that allows us to multiply matrices on each of the modes of the tensor. If $X = [x_{i_1 \dots i_k}] \in \mathbb{F}^{d_1 \times \dots \times d_k}$ and 
\[
A_1 = [a_{ij}^{(1)}] \in \mathbb{F}^{a_1 \times d_1}, \dots, A_k = [a_{ij}^{(k)}] \in \mathbb{F}^{a_k \times d_k},
\]
then $Y= (A_1, \dots, A_k)\cdot X = [y_{i_1 \dots i_k}] \in \mathbb{F}^{a_1 \times \dots \times a_k}$ is the new tensor defined by
\[
y_{i_1 \dots i_k} = \sum_{i_1, \dots, i_k = 1}^{d_1, \dots, d_k} a_{i_1 j_1} \dots a_{i_k j_k} x_{j_1 \dots j_k}.
\]
That is, the matrices $A_i$ act by changes of basis along the $k$ modes of the tensor. Since we are studying symmetric tensors, we impose the condition that $A_1 = \dots = A_k$, otherwise, multilinear matrix multiplication can transform a symmetric tensor $X$ into a non-symmetric tensor $Y$.

The set of $n \times n$ invertible matrices with entries in $\mathbb{F}$ forms a group under matrix multiplication called the general linear group, denoted by GL$_n(\mathbb{F})$. 

\begin{definition}
The group $GL_{d_1, \dots, d_k}(\mathbb{F}) := GL_{d_1}(\mathbb{F}) \times \dots \times GL_{d_k}(\mathbb{F})$ acts on an order-$k$ tensor by changes of basis along the $k$ modes.
\end{definition}

We consider the action of the symmetry group GL$_{2, \dots, 2}(\mathbb{F}_p)$, which does not change the rank of a tensor \cite{deSilvaLim}. Moreover, the action of this group decomposes the set of order-$k$ symmetric tensors into a disjoint union of orbits, where the tensors in each orbit are equivalent under the group action. The orbit of $X$ is the set 
\[
\mathcal{O}_X := \{ (g, g, g) \cdot X \mid g \in \text{GL}_2(\mathbb{F}_p)\}.
\] 
We define the {\bf canonical form} of $X$ to be the minimal element in its orbit with respect to the lexical ordering.

Every order-$k$ symmetric tensor of dimension $n$ may be uniquely associated with a homogeneous polynomial of degree $k$ in $n$ variables \cite{Comon1}. Homogeneous polynomials are also called quantics. Then the problem we consider is equivalent to the problem of determining the minimal number of $p$th powers of linear terms: the Big Waring Problem. This has been considered by many authors, including Ehrenborg \cite{Ehrenborg2} \cite{Ehrenborg3}.

\section{Symmetric Tensors of Format $2 \times 2 \times 2$}

\subsection{$\mathbb{R}$ and $\mathbb{C}$}

Gurevich \cite{Gurevich} and Weinberg \cite{Weinberg} computed the canonical forms of symmetric $2 \times 2 \times 2$ tensors over $\mathbb{R}$ and $\mathbb{C}$ by considering them as homogeneous polynomials. Over $\mathbb{R}$ the canonical forms of $X \in \mathcal{S}^3(\mathbb{R}^2)$ and the symmetric rank are listed in Table \ref{table222symmetricR}.

  \begin{table}
  \[
  \begin{array}{cccc}
  \text{symmetric rank} & \text{canonical form} 
    \\
  \toprule
  0 &    \left[ \begin{array}{cc|cc} 0 & 0 & 0 & 0 \\ 0 & 0 & 0 & 0 \end{array} \right] 
  \\
  \midrule
  1 &   \left[ \begin{array}{cc|cc} 1 & 0 & 0 & 0 \\ 0 & 0 & 0 & 0 \end{array} \right] 
  \\
  \midrule
  2 &   \left[ \begin{array}{cc|cc} 1 & 0 & 0 & 0 \\ 0 & 0 & 0 & 1 \end{array} \right]
  \\
  \midrule
  3 &   \left[ \begin{array}{cc|cc} 0 & 1 & 1 & 0 \\ 1 & 0 & 0 & 0 \end{array} \right]  
  \\[8pt]
  3 &   \left[ \begin{array}{cc|cc} 0 & 1 & 1 & -1 \\ 1 & -1 & -1 & 0 \end{array} \right] 
  \\
  \bottomrule
  \end{array}
  \]
  \medskip
  \caption{Canonical forms of $2 \times 2 \times 2$ symmetric tensors over $\mathbb{R}$}
  \label{table222symmetricR}
  \end{table}

Over $\mathbb{C}$, the list is the same except that we omit the tensor
\[
\left[ \begin{array}{cc|cc} 0 & 1 & 1 & -1 \\ 1 & -1 & -1 & 0 \end{array} \right],
\]
so that there is only one orbit in symmetric rank 3. The maximum symmetric rank is 3 over $\mathbb{R}$ and $\mathbb{C}$. We will see shortly that the maximum rank is larger over some of the prime fields $\mathbb{F}_p$.

\subsection{Algorithms} 
 
The set of $2 \times 2 \times 2$ tensors over the prime field $\mathbb{F}_p$ contains $p^8$ tensors: $p$ many elements in each of the eight entries of the tensor. If we restrict ourselves to the symmetric case, there are only $p^4$ tensors. The algorithm has two main components: generate all the tensors in each rank; apply the group action so that the set of tensors in each rank is decomposed into a disjoint union of orbits. We've used this algorithm in previous work to determine the canonical forms of $2 \times 2 \times 2$ and $2 \times 2 \times 2 \times 2$ tensors over finite fields \cite{BremnerStavrou} . There are two  modification: in the first part, the set of rank-1 tensors must contain only symmetric simple tensors; in the second part, the same group element must act on every mode of the tensor whose orbit we wish to compute (otherwise the orbit will include non-symmetric tensors).

\begin{algorithm} \cite{BremnerStavrou}\\
$\bullet$ The only rank-0 tensor is the zero tensor. The set of symmetric rank-1 tensors is $\{a \otimes a \otimes a \mid a \in \mathbb{F}_p^2\}$ where $a$ is a non-zero vector. 
\\
$\bullet$ Assume that the symmetric rank$_S$-$r$ tensors have already been computed. To compute the symmetric tensors of rank$_S$-$(r+1)$, we consider all sums $X + Y$ where rank$_S(X) = r$ and rank$_S(Y) = 1$. Clearly, rank$_S(X+Y) \leq r+1$, but it is possible that rank$_S(X+Y) < r+1$, so we only retain those $X + Y$ such that rank$_S(X+Y) = r+1$
\\
$\bullet$ Sort the tensors within each symmetric rank in lexical order. 
\\
$\bullet$ For each symmetric rank, perform the following iteration until there are no more tensors in the current symmetric rank:
\\
\indent -- Choose the minimal element of the set of tensors of the current symmetric rank.
\\
\indent -- Compute the orbit of this element under the action of the symmetry group.
\\
\indent -- Remove the elements of this orbit from the set.
\end{algorithm}

See Table \ref{basicalgorithmtable} and Table \ref{groupalgorithmtable} for the pseudocode from \cite{BremnerStavrou} that has been modified for symmetric tensors.

\begin{table}
\begin{itemize}
\item[]
\texttt{flatten}$( x )$
\item[] \qquad
\texttt{return}$( [ \, x_{111}, \, x_{112}, \, x_{121}, \, x_{122}, \, x_{211}, \, x_{212}, \, x_{221}, \, x_{222} \, ] )$
\end{itemize}
\medskip
\begin{itemize}
\item[]
\texttt{outerproduct}$( a, b, c )$
\item[] \qquad
for $i = 1, 2$ do for $j = 1, 2$ do for $k = 1, 2$ do:
set $x_{ijk} \leftarrow a_i a_j a_k$ mod $p$
\item[] \qquad
\texttt{return}( $x$ )
\end{itemize}
\medskip
\begin{itemize}
\item
set $\texttt{vectors} \leftarrow \{ \, [i,j] \mid 0 \le i, j \le p{-}1 \} \setminus \{ [0,0] \}$
\item
set \texttt{arrayset}$[0] \leftarrow \{ [0,0,0,0,0,0,0,0] \}$
\item
set \texttt{arrayset}$[1] \leftarrow \{ \, \}$
\item
for $a$ in \texttt{vectors} do 
  \begin{itemize}
  \item
  set $x \leftarrow \texttt{flatten}( \texttt{outerproduct}( a, a, a ) )$
  \item
  if $x \notin \texttt{arrayset}[0]$ and $x \notin \texttt{arrayset}[1]$ then
  \item[] \quad
  set $\texttt{arrayset}[1] \leftarrow \texttt{arrayset}[1] \cup \{ x \}$
  \end{itemize}
\item
set $r \leftarrow 1$
\item
while $\texttt{arrayset}[r] \ne \{ \, \}$ do:
  \begin{itemize}
  \item
  set $\texttt{arrayset}[r{+}1] \leftarrow \{ \, \}$
  \item
  for $x \in \texttt{arrayset}[r]$ do for $y \in \texttt{arrayset}[1]$ do
    \begin{itemize}
    \item
    set $z \leftarrow [ \, x_1{+}y_1 \, \mathrm{mod}\,p, \, \, \dots, \, x_8{+}y_8 \, \mathrm{mod}\,p \, ]$
    \item
    if $z \notin \texttt{arrayset}[s]$ for $s = 0, \dots, r{+}1$ then
    \item[] \quad
      set $\texttt{arrayset}[r{+}1] \leftarrow \texttt{arrayset}[r{+}1] \cup \{ z \}$
    \end{itemize}
    \item
    set $r \leftarrow r + 1$
  \end{itemize}
\item
set $\texttt{maximumrank} \leftarrow r - 1$
\end{itemize}
\bigskip
\caption{Algorithm to generate tensors in each rank (pseudocode)}
\label{basicalgorithmtable}
\end{table}

\begin{table}
\begin{itemize}
\item[]
\texttt{unflatten}$( x )$
\item[] \quad
set $t \leftarrow 0$
\item[] \quad
for $i = 1, 2$ do for $j = 1, 2$ do for $k = 1, 2$ do:
set $t \leftarrow t + 1$;
set $y_{ijk} \leftarrow x_t$
\item[] \quad
$\texttt{return}( y )$
\end{itemize}
\medskip
\begin{itemize}
\item[]
\texttt{groupaction}$( g, x, m )$
\item[] \quad
set $y \leftarrow \texttt{unflatten}( x )$
\item[] \quad
if $m = 1$ then
for $j = 1,2$ do for $k = 1,2$ do:
\item[] \quad \quad
set $v \leftarrow [ \, y_{1jk}, \, y_{2jk} \, ]$;
set $w \leftarrow [ \, g_{11} v_1 {+} g_{12} v_2 \, \mathrm{mod}\,p, \, g_{21} v_1 {+} g_{22} v_2 \, \mathrm{mod}\,p \, ]$
\item[] \quad \quad
for $i = 1, 2$ do: set $y_{ijk} \leftarrow w_i$
\item[] \quad
if $m = 2$ then \dots \emph{(similar for second subscript)}
\item[] \quad
if $m = 3$ then \dots \emph{(similar for third subscript)}
\item[] \quad
\texttt{return}( \texttt{flatten}( $y$ ) )
\end{itemize}
\medskip
\begin{itemize}
\item[]
\texttt{smallorbit}$( x )$
\item[] \quad
set $\texttt{result} \leftarrow \{\,\}$
\item[] \quad
for $a \in GL_2(\mathbb{F}_p)$ do:
\item[] \quad \quad
set $y \leftarrow \texttt{groupaction}( a, x, 1 )$
\item[] \quad \quad
set $z \leftarrow \texttt{groupaction}( a, y, 2 )$
\item[] \quad \quad 
set $w \leftarrow \texttt{groupaction}( a, z, 3 )$
\item[] \quad
set $\texttt{result} \leftarrow \texttt{result} \cup \{ w \}$
\item[] \quad
\texttt{return}( \texttt{result} )
\end{itemize}
\medskip
\begin{itemize}
\item
for $r = 0,\dots,\texttt{maximumrank}$ do:
\item[] \quad
set $\texttt{representatives}[r] \leftarrow \{\,\}$;
set $\texttt{remaining} \leftarrow \texttt{arrayset}[r]$
\item[] \quad
while $\texttt{remaining} \ne \{\,\}$ do:
\item[] \quad \quad
set $x \leftarrow \texttt{remaining}[1]$;
set $\texttt{xorbit} \leftarrow \texttt{largeorbit}( x )$
\item[] \quad \quad
append $\texttt{xorbit}[1]$ to $\texttt{representatives}[r]$
\item[] \quad \quad
set $\texttt{remaining} \leftarrow \texttt{remaining} \setminus \texttt{xorbit}$
\end{itemize}
\medskip
\caption{Algorithm for group action (pseudocode)}
\label{groupalgorithmtable}
\end{table}

\subsection{$\mathbb{F}_2$}

Over $\mathbb{F}_2$, there are 256 tensors, 16 of which are symmetric. There are three symmetric simple tensors,
\[
\left[
\begin{array}{cc|cc}
1 & 0 & 0 & 0 \\
0 & 0 & 0 & 0
\end{array}
\right],
\quad
\left[
\begin{array}{cc|cc}
0 & 0 & 0 & 0 \\
0 & 0 & 0 & 1
\end{array}
\right],
\quad
\left[
\begin{array}{cc|cc}
1 & 1 & 1 & 1 \\
1 & 1 & 1 & 1
\end{array}
\right].
\]
We emphasize that not all 16 symmetric tensors can be written as a sum of symmetric simple tensors over $\mathbb{F}_2$. The following 8 tensors do not have a symmetric decomposition: 
\begin{align*}
&
\left[
\begin{array}{cc|cc}
0 & 1 & 1 & 0 \\
1 & 0 & 0 & 0 
\end{array}
\right],
\quad
\left[
\begin{array}{cc|cc}
1 & 1 & 1 & 0 \\
1 & 0 & 0 & 0 
\end{array}
\right],
\quad
\left[
\begin{array}{cc|cc}
0 & 1 & 1 & 0 \\
1 & 0 & 0 & 1 
\end{array}
\right],
\quad
\left[
\begin{array}{cc|cc}
1 & 1 & 1 & 0 \\
1 & 0 & 0 & 1 
\end{array}
\right],
\\
&
\left[
\begin{array}{cc|cc}
0 & 0 & 0 & 1 \\
0 & 1 & 1 & 0 
\end{array}
\right],
\quad
\left[
\begin{array}{cc|cc}
1 & 0 & 0 & 1 \\
0 & 1 & 1 & 0 
\end{array}
\right],
\quad
\left[
\begin{array}{cc|cc}
0 & 0 & 0 & 1 \\
0 & 1 & 1 & 1 
\end{array}
\right],
\quad
\left[
\begin{array}{cc|cc}
1 & 0 & 0 & 1 \\
0 & 1 & 1 & 1 
\end{array}
\right].
\end{align*}

The maximum symmetric rank over $\mathbb{F}_2$ is 3. The number of symmetric tensors in each symmetric rank and the approximate percentages are listed below. The percentages add to 50\%, since only half of the tensors have a decomposition into a sum of symmetric simple tensors.
\[
\begin{array}{lrrrr}
\text{rank} & 0 & 1 & 2 & 3 \\
\text{number} & 1 & 3 & 3 & 1 \\
\text{$\approx$ $\%$} & 6.25\% & 18.75\% & 18.75\% & 6.25\%
\end{array}
\]

The symmetric ranks, orders of each orbit, and the minimal representatives of each orbit are given in Table \ref{table222symmetricmod2}.
 
  \begin{table}
  \[
  \begin{array}{cccc}
  \text{symmetric rank} & \text{orbit size} &  \text{canonical form} 
    \\
  \toprule
  0 & 1 &    \left[ \begin{array}{cc|cc} 0 & 0 & 0 & 0 \\ 0 & 0 & 0 & 0 \end{array} \right] 
  \\
  \midrule
  1 & 3 &   \left[ \begin{array}{cc|cc} 0 & 0 & 0 & 0 \\ 0 & 0 & 0 & 1 \end{array} \right] 
  \\
  \midrule
  2 & 3 &  \left[ \begin{array}{cc|cc} 0 & 1 & 1 & 1 \\ 1 & 1 & 1 & 1 \end{array} \right]
  \\
  \midrule
  3 & 1 &   \left[ \begin{array}{cc|cc} 0 & 1 & 1 & 1 \\ 1 & 1 & 1 & 0 \end{array} \right]  
  \\
  \bottomrule
  \end{array}
  \]
  \medskip
  \caption{Canonical forms of $2 \times 2 \times 2$ symmetric tensors over $\mathbb{F}_2$}
  \label{table222symmetricmod2}
  \end{table}

\subsection{$\mathbb{F}_3$}

There are 6,561 tensors over $\mathbb{F}_3$, 81 of which are symmetric. Unlike in the previous case, all 81 symmetric tensors can be generated as a sum of simple symmetric tensors. The maximum symmetric rank in this case is 4. The number of symmetric tensors in each symmetric rank and the approximate percentages are listed below.
\[
\begin{array}{lrrrrr}
\text{rank} & 0 & 1 & 2 & 3 & 4 \\
\text{number} & 1 & 8 & 24 & 32 & 16 \\
\text{$\approx$ $\%$} & 1.23\% & 9.88\% & 29.63\% & 39.51\% & 19.75\%
\end{array}
\]

The symmetric ranks, orders of each orbit, and the minimal representatives of each orbit are given in Table \ref{table222symmetricmod3}. 

  \begin{table}
  \[
  \begin{array}{cccc}
  \text{symmetric rank} & \text{orbit size} &  \text{canonical form} 
    \\
  \toprule
  0 & 1 &    \left[ \begin{array}{cc|cc} 0 & 0 & 0 & 0 \\ 0 & 0 & 0 & 0 \end{array} \right] 
  \\
  \midrule
  1 & 8 &   \left[ \begin{array}{cc|cc} 0 & 0 & 0 & 0 \\ 0 & 0 & 0 & 1 \end{array} \right] 
  \\
  \midrule
  2 & 24 & \left[ \begin{array}{cc|cc} 0 & 1 & 1 & 0 \\ 1 & 0 & 0 & 1 \end{array} \right]
  \\
  \midrule
  3 & 8 & \left[ \begin{array}{cc|cc} 0 & 0 & 0 & 1 \\ 0 & 1 & 1 & 0 \end{array} \right]  
  \\[7pt]
  3 & 24 &   \left[ \begin{array}{cc|cc} 0 & 1 & 1 & 0 \\ 1 & 0 & 0 & 2 \end{array} \right]
  \\
  \midrule
  4 & 16 &   \left[ \begin{array}{cc|cc} 0 & 0 & 0 & 1 \\ 0 & 1 & 1 & 1 \end{array} \right]
  \\
  \bottomrule
  \end{array}
  \]
  \medskip
  \caption{Canonical forms of $2 \times 2 \times 2$ symmetric tensors over $\mathbb{F}_3$}
  \label{table222symmetricmod3}
  \end{table}

\subsection{$\mathbb{F}_5$}
There are 390,625 tensors over $\mathbb{F}_3$, 625 of which are symmetric. The maximum symmetric rank in this case is 3. The number of tensors in each symmetric rank and the approximate percentages are listed below. 
\[
\begin{array}{lrrrrr}
\text{rank} & 0 & 1 & 2 & 3 \\
\text{number} & 1 & 24 & 240 & 360 \\
\text{$\approx$ $\%$} & 0.16\% & 3.84\% & 38.40\% & 57.60\% 
\end{array}
\]

The symmetric ranks, orders of each orbit, and the minimal representatives of each orbit are given in Table \ref{table222symmetricmod5}. From now on, the canonical forms will be displayed in flattened form.

  \begin{table}
  \[
  \begin{array}{ccc}
  \text{symmetric rank} & \text{orbit size} & \text{canonical form (flattened)} 
    \\
  \toprule
  0 & 1 &    \begin{bmatrix} 0 & 0 & 0 & 0 & 0 & 0 & 0 & 0 \end{bmatrix}  
  \\
  \midrule
  1 & 24 &   \begin{bmatrix} 0 & 0 & 0 & 0 & 0 & 0 & 0 & 1 \end{bmatrix}  
  \\
  \midrule
  2 & 240 & \begin{bmatrix} 0 & 1 & 1 & 0 & 1 & 0 & 0 & 1 \end{bmatrix} 
  \\
  \midrule
  3 & 120 &   \begin{bmatrix} 0 & 0 & 0 & 1 & 0 & 1 & 1 & 0 \end{bmatrix}   
  \\
  3 & 80 &  \begin{bmatrix} 0 & 1 & 1 & 0 & 1 & 0 & 0 & 2 \end{bmatrix} 
  \\
  3 & 160 &   \begin{bmatrix} 1 & 0 & 0 & 1 & 0 & 1 & 1 & 2 \end{bmatrix} 
  \\
  \bottomrule
  \end{array}
  \]
  \medskip
  \caption{Canonical forms of $2 \times 2 \times 2$ symmetric tensors over $\mathbb{F}_5$}
  \label{table222symmetricmod5}
  \end{table}

\subsection{$\mathbb{F}_7$}
There are 5,764,801 tensors over $\mathbb{F}_7$, 2,401 of which are symmetric. The maximum symmetric rank in this case is 5, and we also see that the number of orbits has jumped to 14. The number of symmetric tensors in each symmetric rank and the approximate percentages are listed below.
\[
\begin{array}{lrrrrrr}
\text{rank} & 0 & 1 & 2 & 3 & 4 & 5 \\
\text{number} & 1 & 16 & 128 & 688 & 1232 & 336 \\
\text{$\approx$ $\%$} & 0.04\% & 0.67\% & 5.33\% & 28.65\% & 51.31\% & 13.99\%
\end{array}
\]

The symmetric ranks, orders of each orbit, and the minimal representatives of each orbit are given in Table \ref{table222symmetricmod7}.

  \begin{table}
  \[
  \begin{array}{cccc}
  \text{symmetric rank} & \text{orbit size} &  \text{canonical form (flattened)} 
    \\
  \toprule
  0 & 1 &    \begin{matrix} 0 & 0 & 0 & 0 & 0 & 0 & 0 & 0 \end{matrix}
  \\
  \midrule
  1 & 16 &   \begin{matrix} 0 & 0 & 0 & 0 & 0 & 0 & 0 & 1 \end{matrix} 
  \\
  \midrule
  2 & 16 & \begin{matrix} 0 & 0 & 0 & 0 & 0 & 0 & 0 & 2 \end{matrix}
  \\
  2 & 112 & \begin{matrix} 0 & 1 & 1 & 0 & 1 & 0 & 0 & 2 \end{matrix}
  \\
  \midrule
  3 & 16 & \begin{matrix} 0 & 0 & 0 & 0 & 0 & 0 & 0 & 3 \end{matrix}  
  \\
  3 & 112 &   \begin{matrix} 0 & 1 & 1 & 0 & 1 & 0 & 0 & 1 \end{matrix}
  \\
  3 & 336 &   \begin{matrix} 0 & 1 & 1 & 0 & 1 & 0 & 0 & 3 \end{matrix}
  \\
  3 & 224 &   \begin{matrix} 1 & 0 & 0 & 0 & 0 & 0 & 0 & 2 \end{matrix}
    \\
  \midrule
  4 & 336 & \begin{matrix} 0 & 0 & 0 & 1 & 0 & 1 & 1 & 0 \end{matrix}  
  \\
  4 & 112 &   \begin{matrix} 0 & 1 & 1 & 0 & 1 & 0 & 0 & 6 \end{matrix}
  \\
  4 & 336 &   \begin{matrix} 0 & 1 & 1 & 0 & 1 & 0 & 0 & 6 \end{matrix}
  \\
  4 & 224 &   \begin{matrix} 1 & 0 & 0 & 0 & 0 & 0 & 0 & 3 \end{matrix}
  \\
  4 & 224 &   \begin{matrix} 1 & 0 & 0 & 1 & 0 & 1 & 1 & 2 \end{matrix}
  \\
  \midrule
  5 & 336 &   \begin{matrix} 0 & 1 & 1 & 0 & 1 & 0 & 0 & 5 \end{matrix}
  \\
  \bottomrule
  \end{array}
  \]
  \medskip
  \caption{Canonical forms of $2 \times 2 \times 2$ symmetric tensors over $\mathbb{F}_7$}
  \label{table222symmetricmod7}
  \end{table}

\subsection{$\mathbb{F}_{11}$}
There are 214,358,881 tensors over $\mathbb{F}_{11}$, 14,641 of which are symmetric. The maximum symmetric rank in this case is 3. The number of symmetric tensors in each symmetric rank and the approximate percentages are listed below.
\[
\begin{array}{lrrrr}
\text{rank} & 0 & 1 & 2 & 3 \\
\text{number} & 1 & 120 & 6600 & 7920  \\
\text{$\approx$ $\%$} & 0.01\% & 0.82\% & 45.08\% & 54.09\%
\end{array}
\]

The symmetric ranks, orders of each orbit, and the minimal representatives of each orbit are given in Table \ref{table222symmetricmod11}.

  \begin{table}
  \[
  \begin{array}{cccc}
  \text{symmetric rank} & \text{orbit size} &  \text{canonical form (flattened)} 
    \\
  \toprule
  0 & 1 &    \begin{matrix} 0 & 0 & 0 & 0 & 0 & 0 & 0 & 0 \end{matrix}
  \\
  \midrule
  1 & 120 &   \begin{matrix} 0 & 0 & 0 & 0 & 0 & 0 & 0 & 1 \end{matrix} 
  \\
  \midrule
  2 & 6600 & \begin{matrix} 0 & 1 & 1 & 0 & 1 & 0 & 0 & 1 \end{matrix}
  \\
  \midrule
  3 & 1320 & \begin{matrix} 0 & 0 & 0 & 1 & 0 & 1 & 1 & 0 \end{matrix}  
  \\
  3 & 2200 &   \begin{matrix} 0 & 1 & 1 & 0 & 1 & 0 & 0 & 2 \end{matrix}
  \\
  3 & 4400 &   \begin{matrix} 1 & 0 & 0 & 1 & 0 & 1 & 1 & 2 \end{matrix}
  \\
  \bottomrule
  \end{array}
  \]
  \medskip
  \caption{Canonical forms of $2 \times 2 \times 2$ symmetric tensors over $\mathbb{F}_{11}$}
  \label{table222symmetricmod11}
  \end{table}

\subsection{$\mathbb{F}_{13}$}
There are 815,730,721 tensors over $\mathbb{F}_{13}$, 28,561 of which are symmetric.  The maximum symmetric rank in this case is 4, and there are 14 orbits. The number of symmetric tensors in each symmetric rank and the approximate percentages are listed below.
\[
\begin{array}{lrrrrrr}
\text{rank} & 0 & 1 & 2 & 3 & 4 \\
\text{number} & 1 & 56 & 1568 & 16016 & 10920 \\
\text{$\approx$ $\%$} & 0.00\% & 0.20\% & 5.49\% & 56.08\% & 38.23\%
\end{array}
\]

The symmetric ranks, orders of each orbit, and the minimal representatives of each orbit are given in Table \ref{table222symmetricmod13}. 

  \begin{table}
  \[
  \begin{array}{cccc}
  \text{symmetric rank} & \text{orbit size} &  \text{canonical form (flattened)} 
    \\
  \toprule
  0 & 1 &    \begin{matrix} 0 & 0 & 0 & 0 & 0 & 0 & 0 & 0 \end{matrix}
  \\
  \midrule
  1 & 56 &   \begin{matrix} 0 & 0 & 0 & 0 & 0 & 0 & 0 & 1 \end{matrix} 
  \\
  \midrule
  2 & 56 & \begin{matrix} 0 & 0 & 0 & 0 & 0 & 0 & 0 & 2 \end{matrix}
  \\
  2 & 56 & \begin{matrix} 0 & 0 & 0 & 0 & 0 & 0 & 0 & 4 \end{matrix}
  \\
  2 &14 56 & \begin{matrix} 0 & 1 & 1 & 0 & 1 & 0 & 0 & 3 \end{matrix}
  \\
  \midrule
  3 & 4368 & \begin{matrix} 0 & 1 & 1 & 0 & 1 & 0 & 0 & 2 \end{matrix}  
  \\
  3 & 1456 &   \begin{matrix} 0 & 1 & 1 & 0 & 1 & 0 & 0 & 4 \end{matrix}
  \\
  3 & 4368 &   \begin{matrix} 0 & 1 & 1 & 0 & 1 & 0 & 0 & 5 \end{matrix}
  \\
  3 & 2912 &   \begin{matrix} 1 & 0 & 0 & 0 & 0 & 0 & 0 & 2 \end{matrix}
   \\
  3 & 2912 &   \begin{matrix} 1 & 0 & 0 & 0 & 0 & 0 & 0 & 4 \end{matrix}
   \\
  \midrule
  4 & 2184 & \begin{matrix} 0 & 0 & 0 & 1 & 0 & 1 & 1 & 0 \end{matrix}  
  \\
  4 & 1456 &   \begin{matrix} 0 & 1 & 1 & 0 & 1 & 0 & 0 & 1 \end{matrix}
  \\
  4 & 4368 &   \begin{matrix} 0 & 1 & 1 & 0 & 1 & 0 & 0 & 6 \end{matrix}
  \\
  4 & 2912 &   \begin{matrix} 1 & 0 & 0 & 1 & 0 & 1 & 1 & 6 \end{matrix}
  \\
  \bottomrule
  \end{array}
  \]
  \medskip
  \caption{Canonical forms of $2 \times 2 \times 2$ symmetric tensors over $\mathbb{F}_{13}$}
  \label{table222symmetricmod13}
  \end{table}

\subsection{$\mathbb{F}_{17}$}
There are 6,975,757,441 tensors over $\mathbb{F}_{13}$, 83,521 of which are symmetric.  The maximum symmetric rank in this case is 3, and there are 6 orbits. The number of symmetric tensors in each symmetric rank and the approximate percentages are listed below.
\[
\begin{array}{lrrrrrr}
\text{rank} & 0 & 1 & 2 & 3 \\
\text{number} & 1 & 288 & 39168 & 44064 \\
\text{$\approx$ $\%$} & 0.00\% & 0.34\% & 46.90\% & 52.76\%
\end{array}
\]

The symmetric ranks, orders of each orbit, and the minimal representatives of each orbit are given in Table \ref{table222symmetricmod17}.

  \begin{table}
  \[
  \begin{array}{cccc}
  \text{symmetric rank} & \text{orbit size} &  \text{canonical form (flattened)} 
    \\
  \toprule
  0 & 1 &    \begin{matrix} 0 & 0 & 0 & 0 & 0 & 0 & 0 & 0 \end{matrix}
  \\
  \midrule
  1 & 56 &   \begin{matrix} 0 & 0 & 0 & 0 & 0 & 0 & 0 & 1 \end{matrix} 
  \\
  \midrule
  2 & 56 & \begin{matrix} 0 & 1 & 1 & 0 & 1 & 0 & 0 & 1 \end{matrix}
  \\
  \midrule
  3 & 4368 & \begin{matrix} 0 & 0 & 0 & 1 & 0 & 1 & 1 & 0 \end{matrix}  
  \\
  3 & 1456 &   \begin{matrix} 0 & 1 & 1 & 0 & 1 & 0 & 0 & 3 \end{matrix}
  \\
  3 & 4368 &   \begin{matrix} 1 & 0 & 0 & 1 & 0 & 1 & 1 & 1 \end{matrix}
  \\
  \bottomrule
  \end{array}
  \]
  \medskip
  \caption{Canonical forms of $2 \times 2 \times 2$ symmetric tensors over $\mathbb{F}_{17}$}
  \label{table222symmetricmod17}
  \end{table}

For $p = 19$ there are almost 17 billion tensors, where over 130,000 of these are symmetric. For this and larger primes, our program does not have sufficient memory.


\section{Symmetric Tensors of Format $2 \times 2 \times 2 \times 2$}

\subsection{$\mathbb{R}$ and $\mathbb{C}$}
Weinberg \cite{Weinberg} lists the canonical forms of $2 \times 2 \times 2 \times 2$ symmetric tensors based on the discussion from Gurevich \cite{Gurevich}. Recall we represent these fourth order tensors as $X = [x_{ijkl}]$ where $i,j,k,l = 1, 2$. Over $\mathbb{R}$ there are 9 canonical forms:
\begin{itemize}
\item[1.]{$x_{1111} = 1, \ x_{\sigma(1122)} = \mu, \ x_{2222} = 1, \ \mu < -\frac{1}{3}, \ x_{ijkl} = 0$ otherwise;}
\item[2.]{$x_{1111} = \alpha, \ x_{\sigma(1122)} = \alpha \mu, \ x_{2222} = \alpha, \ \mu > -\frac{1}{3}, \ \alpha = \pm1, \ x_{ijkl} = 0$ otherwise;}
\item[3.]{$x_{1111} = 1, \ x_{\sigma(1122)} = \mu, \ x_{2222} = -1, \ x_{ijkl} = 0$ otherwise; }
\item[4.]{ $x_{\sigma(1112)} = \alpha, \ x_{ijkl} = 0$ otherwise;}
\item[5.]{ $x_{ijkl} = 0$;}
\item[6.]{ $x_{\sigma(1122)} = \pm 1, \ x_{2222} = \pm 1, \ x_{ijkl} = 0$ otherwise;}
\item[7.]{ $x_{\sigma(1122)} = \alpha, \ x_{ijkl} = 0$ otherwise;}
\item[8.]{ $x_{1111} = \gamma, \ x_{\sigma(1122)} = \frac{\gamma}{3}, \ x_{2222} = \gamma, \gamma = 1, \ x_{ijkl} = 0$ otherwise;}
\item[9.]{$x_{1111} = \pm 1, \ x_{ijkl} = 0$ otherwise.}
\end{itemize}

Here, $\sigma \in S_4$ denotes any permutation of the indicated subscripts. If we consider the group action GL$_{2,2,2,2}(\mathbb{C})$ then we get 6 canonical forms:
\begin{itemize}
\item[1.]{$x_{1111} = 1, \ x_{\sigma(1122)} = 6\mu, \ x_{2222} = 1, \ \mu \ne \pm \frac{1}{3}, \ x_{ijkl} = 0$ otherwise; }
\item[2.]{$x_{\sigma(1122)} = 6, \ x_{2222} = 1, \ x_{ijkl} = 0$ otherwise;}
\item[3.]{$x_{\sigma(1122)} = 6, \ x_{ijkl} = 0$ otherwise;}
\item[4.]{$x_{\sigma(1112)} = 4, \ x_{ijkl} = 0$ otherwise;}
\item[5.]{$x_{1111} = 1, \ x_{ijkl} = 0$ otherwise;}
\item[6.]{$x_{ijkl} = 0$.}
\end{itemize}

The algorithms we used for the $2 \times 2 \times 2$ case generalize to the $2 \times 2 \times 2 \times 2$ case by making a few simple modifications. In particular, we consider the outer product of four non-zero vectors (instead of 3), and we add an extra index $\ell \in \{1,2\}$ for the group action on the fourth mode.

\subsection{$\mathbb{F}_{2}$}
There are 65,536 tensors over $\mathbb{F}_2$, 32 of which are symmetric. Just as we experienced for the $2 \times 2 \times 2$ case over $\mathbb{F}_2$, not all the symmetric tensors can be represented as a sum of symmetric simple tensors. In this case, only 8 out of 32 symmetric tensors can be represented as a sum of symmetric simple tensors. The maximum symmetric rank in this case is 3, and there are 4 orbits. The number of symmetric tensors in each symmetric rank and the approximate percentages are listed below.
\[
\begin{array}{lrrrrrr}
\text{rank} & 0 & 1 & 2 & 3 \\
\text{number} & 1 & 3 & 3 & 1 \\
\text{$\approx$ $\%$} & 3.125\% & 9.375\% & 9.375\% & 3.125\%
\end{array}
\]

The symmetric ranks, orders of each orbit, and the minimal representatives of each orbit are given in Table \ref{table2222symmetricmod2}. To improve legibility, we use a period (.) in place of zero (0).

  \begin{table}
  \[
  \begin{array}{ccccccccccccccccccc}
  \text{sym rank} & \text{orbit size} &\multicolumn{16}{c}{\text{canonical form (flattened)}} 
    \\
  \toprule
 0 & 1 & . & . & . & . & . & . & . & . & . & . & . & . & . & . & . & .  
  \\
  \midrule
  1 & 56 &  . & . & . & . & . & . & . & . & . & . & . & . & . & . & . & 1  
  \\
  \midrule
  2 & 56 &  . & 1 & . & . & . & . & . & . & . & . & . & . & . & . & . & 1 
  \\
  \midrule
  3 & 4368 &  . & 1 & . & . & . & . & . & . & . & . & . & . & . & . & 1 & .   
  \\
  \bottomrule
  \end{array}
  \]
  \medskip
  \caption{Canonical forms of $2 \times 2 \times 2 \times 2$ symmetric tensors over $\mathbb{F}_2$}
  \label{table2222symmetricmod2}
  \end{table}

\subsection{$\mathbb{F}_3$}
There are 43,046,721 tensors over $\mathbb{F}_3$, 243 of which are symmetric. Only 81 out of the 243 symmetric tensors can be written as a sum of symmetric simple tensors. The maximum symmetric rank in this case is 8, and there are 15 orbits. The number of symmetric tensors in each symmetric rank and the approximate percentages are listed below.
\[
\begin{array}{lrrrrrrrrrrr}
\text{rank} & 0 & 1 & 2 & 3 & 4 & 5 & 6 & 7 & 8 \\
\text{number} & 1 & 4 & 10 & 16 & 19 & 16 & 10 & 4 & 1 \\
\text{$\approx$ $\%$} & 0.41\% & 1.65\% & 4.12\% & 6.58\% & 7.82\% & 6.58\% & 4.12\% &  1.65\% & 0.41\%
\end{array}
\]

The symmetric ranks, orders of each orbit, and the minimal representatives of each orbit are given in Table \ref{table2222symmetricmod3}.

  \begin{table}
  \[
  \begin{array}{ccccccccccccccccccc}
  \text{sym rank} & \text{orbit size} &\multicolumn{16}{c}{\text{canonical form (flattened)}} 
    \\
  \toprule
 0 & 1 & . & . & . & . & . & . & . & . & . & . & . & . & . & . & . & .  
  \\
  \midrule
  1 & 4 &  . & . & . & . & . & . & . & . & . & . & . & . & . & . & . & 1  
  \\
  \midrule
  2 & 4 &  . & . & . & . & . & . & . & . & . & . & . & . & . & . & . & 2 
  \\
  2 & 6 &  1 & . & . & . & . & . & . & . & . & . & . & . & . & . & . & 1 
  \\
  \midrule
  3 & 4 &  . & . & . & 2 & . & 2 & 2 & . & . & 2 & 2 & . & 2 & . & . & 2   
  \\
  3 & 12 &  . & 1 & 1 & . & 1 & . & . & 1 & 1 & . & . & 1 & . & 1 & 1 & .   
  \\
  \midrule
  4 & 1 & . & . & . & 2 & . & 2 & 2 & . & . & 2 & 2 & . & 2 & . & . & . 
  \\
  4 & 12 & . & 1 & 1 & . & 1 & . & . & 1 & 1 & . & . & 1 & . & 1 & 1 & 1
  \\
  4 & 6 & 1 & . & . & 1 & . & 1 & 1 & . & . & 1 & 1 & . & 1 & . & . & 1
  \\
  \midrule
  5 & 4 & . & . & . & 2 & . & 2 & 2 & . & . & 2 & 2 & . & 2 & . & . & 1
  \\
  5 & 12 & . & 1 & 1 & . & 1 & . & . & 1 & 1 & . & . & 1 & . & 1 & 1 & 2
  \\
  \midrule
  6 & 4 & . & . & . & 1 & . & 1 & 1 & . & . & 1 & 1 & . & 1 & . & . & 1
  \\
  6 & 6 & 1 & . & . & 2 & . & 2 & 2 & . & . & 2 & 2 & . & 2 & . & . & 1
  \\
  \midrule
  7 & 4 & . & . & . & 1 & . & 1 & 1 & . & . & 1 & 1 & . & 1 & . & . & 2
  \\
  \midrule
  8 & 1 & . & . & . & 1 & . & 1 & 1 & . & . & 1 & 1 & . & 1 & . & . & .
  \\
  \bottomrule
  \end{array}
  \]
  \medskip
  \caption{Canonical forms of $2 \times 2 \times 2 \times 2$ symmetric tensors over $\mathbb{F}_3$}
  \label{table2222symmetricmod3}
  \end{table}

\subsection{$\mathbb{F}_5$}
There are 152,587,890,625 tensors over $\mathbb{F}_5$, 3,125 of which are symmetric. In this case, all symmetric tensors have a symmetric decomposition. The maximum symmetric rank in this case is 10, and there are 52 orbits. The number of symmetric tensors in each symmetric rank and the approximate percentages are listed below.
\[
\begin{array}{lrrrrrr}
\text{rank} & 0 & 1 & 2 & 3 & 4 & 5 \\
\text{number} & 1 & 6 & 21 & 56 & 126 & 240 \\
\text{$\approx$ $\%$} & 0.03\% & 0.19\% & 0.67\% & 1.79\% & 4.03\% & 7.68\% \\ \hline
\text{rank} & 6 & 7 & 8 & 9 & 10 & \\
\text{number} & 395 & 570 & 690 & 660 & 360 &  \\
\text{$\approx$ $\%$} & 12.64\% & 18.24\% & 22.08\% & 21.12\% & 11.52\% & 
\end{array}
\]

The symmetric ranks, orders of each orbit, and the minimal representatives of each orbit are given in Tables \ref{table2222symmetricmod5-1} and \ref{table2222symmetricmod5-2}. 

  \begin{table}
  \[
  \begin{array}{ccccccccccccccccccc}
  \text{sym rank} & \text{orbit size} &\multicolumn{16}{c}{\text{canonical form (flattened)}} 
    \\
  \toprule
 0 & 1 & . & . & . & . & . & . & . & . & . & . & . & . & . & . & . & .  
  \\
  \midrule
  1 & 6 &  . & . & . & . & . & . & . & . & . & . & . & . & . & . & . & 1  
  \\
  \midrule
  2 & 6 &  . & . & . & . & . & . & . & . & . & . & . & . & . & . & . & 2 
  \\
  2 & 15 &  1 & . & . & . & . & . & . & . & . & . & . & . & . & . & . & 1 
  \\
  \midrule
  3 & 6 &  . & . & . & . & . & . & . & . & . & . & . & . & . & . & . & 3  
  \\
  3 & 39 &  1 & . & . & . & . & . & . & . & . & . & . & . & . & . & . & 2   
  \\
   3 & 20 &  2 & . & . & 2 & . & 2 & 2 & . & . & 2 & 2 & . & 2 & . & . & 3   
  \\
  \midrule
  4 & 6 & . & . & . & . & . & . & . & . & . & . & . & . & . & . & . & 4
  \\
  4 & 30 & 1 & . & . & . & . & . & . & . & . & . & . & . & . & . & . & 3
  \\
  4 & 15 & 2 & . & . & . & . & . & . & . & . & . & . & . & . & . & . & 2
  \\
   4 & 60 & 2 & . & . & 2 & . & 2 & 2 & . & . & 2 & 2 & . & 2 & . & . & 4
  \\
   4 & 15 & 3 & . & . & 2 & . & 2 & 2 & . & . & 2 & 2 & . & 2 & . & . & 3
  \\
  \midrule
  5 & 60 & . & . & . & 1 & . & 1 & 1 & . & . & 1 & 1 & . & 1 & . & . & 4
  \\
  5 & 60 & . & . & . & 2 & . & 2 & 2 & . & . & 2 & 2 & . & 2 & . & . & 2
  \\
  5 & 30 & . & 1 & 1 & . & 1 & . & . & 1 & 1 & . & . & 1 & . & 1 & 1 & .
  \\
  5 & 30 & . & 1 & 1 & . & 1 & . & . & 4 & 1 & . & . & 4 & . & 4 & 4 & .
  \\
  5 & 60 & . & 1 & 1 & . & 1 & . & . & 4 & 1 & . & . & 4 & . & 4 & 4 & 4
  \\
  \midrule
  6 & 30 & . & . & . & 1 & . & 1 & 1 & . & . & 1 & 1 & . & 1 & . & . & .
  \\
  6 & 60 & . & . & . & 2 & . & 2 & 2 & . & . & 2 & 2 & . & 2 & . & . & 3  
  \\
   6 & 120 & . & 1 & 1 & . & 1 & . & . & 1 & 1 & . & . & 1 & . & 1 & 1 & 1
  \\
  6 & 60 & . & 1 & 1 & . & 1 & . & . & 3 & 1 & . & . & 3 & . & 3 & 3 & .
  \\
   6 & 60 & . & 1 & 1 & . & 1 & . & . & 4 & 1 & . & . & 4 & . & 4 & 4 & 1
  \\
  6 & 15 & 1 & . & . & 1 & . & 1 & 1 & . & . & 1 & 1 & . & 1 & . & . & 1 
  \\
   6 & 20 & 1 & . & . & 1 & . & 1 & 1 & . & . & 1 & 1 & . & 1 & . & . & 4 
  \\
  6 & 30 & 1 & . & . & 2 & . & 2 & 2 & 2 & . & 2 & 2 & 2 & 2 & 2 & 2 & 1
  \\
  \midrule
  7 & 60 & . & . & . & 1 & . & 1 & 1 & . & . & 1 & 1 & . & 1 & . & . & 1
  \\
  7 & 60 & . & . & . & 2 & . & 2 & 2 & . & . & 2 & 2 & . & 2 & . & . & 4
  \\
  7 & 60 & . & 1 & 1 & . & 1 & . & . & 1 & 1 & . & . & 1 & . & 1 & 1 & 2
  \\
  7 & 120 & . & 1 & 1 & . & 1 & . & . & 2 & 1 & . & . & 2 & . & 2 & 2 & 4
  \\
  7 & 120 & . & 1 & 1 & . & 1 & . & . & 4 & 1 & . & . & 4 & . & 4 & 4 & 2
  \\
  7 & 60 & 1 & . & . & 1 & . & 1 & 1 & . & . & 1 & 1 & . & 1 & . & . & 2 
  \\
  7 & 60 & 1 & . & . & 1 & . & 1 & 1 & 2 & . & 1 & 1 & 2 & 1 & 2 & 2 & 2
  \\
  7 & 30 & 2 & . & . & 1 & . & 1 & 1 & 2 & . & 1 & 1 & 2 & 1 & 2 & 2 & 2
    \\
  \bottomrule
  \end{array}
  \]
  \medskip
  \caption{Canonical forms of $2 \times 2 \times 2 \times 2$ symmetric tensors over $\mathbb{F}_5$}
  \label{table2222symmetricmod5-1}
  \end{table}

   \begin{table}
  \[
  \begin{array}{ccccccccccccccccccc}
  \text{sym rank} & \text{orbit size} &\multicolumn{16}{c}{\text{canonical form (flattened)}} 
    \\
  \toprule
  8 & 60 & . & . & . & 1 & . & 1 & 1 & . & . & 1 & 1 & . & 1 & . & . & 2
  \\
   8 & 30 & . & . & . & 2 & . & 2 & 2 & . & . & 2 & 2 & . & 2 & . & . & 0
  \\
   8 & 60 & . & 1 & 1 & . & 1 & . & . & 1 & 1 & . & . & 1 & . & 1 & 1 & 3 
  \\
   8 & 60 & . & 1 & 1 & . & 1 & . & . & 2 & 1 & . & . & 2 & . & 2 & 2 & 0
  \\
   8 & 120 & . & 1 & 1 & . & 1 & . & . & 3 & 1 & . & . & 3 & . & 3 & 3 & 2
  \\
   8 & 120 & . & 1 & 1 & . & 1 & . & . & 4 & 1 & . & . & 4 & . & 4 & 4 & 3
  \\
   8 & 120 & 1 & . & . & 1 & . & 1 & 1 & 2 & . & 1 & 1 & 2 & 1 & 2 & 2 & 3
  \\
   8 & 60 & 1 & . & . & 2 & . & 2 & 2 & . & . & 2 & 2 & . & 2 & . & . & 4 
  \\
   8 & 60 & 1 & . & . & 2 & . & 2 & 2 & 2 & . & 2 & 2 & 2 & 2 & 2 & 2 & 3
  \\  \midrule
  9 & 60 & . & . & . & 1 & . & 1 & 1 & . & . & 1 & 1 & . & 1 & . & . & 3
  \\
   9 & 60 & . & . & . & 2 & . & 2 & 2 & . & . & 2 & 2 & . & 2 & . & . & 1 
  \\
   9 & 120 & . & 1 & 1 & . & 1 & . & . & 1 & 1 & . & . & 1 & . & 1 & 1 & 4 
  \\
   9 & 120 & . & 1 & 1 & . & 1 & . & . & 2 & 1 & . & . & 2 & . & 2 & 2 & 1
  \\
   9 & 120 & . & 1 & 1 & . & 1 & . & . & 3 & 1 & . & . & 3 & . & 3 & 3 & 3
  \\
   9 & 120 & 1 & . & . & . & . & . & . & 1 & . & . & . & 1 & . & 1 & 1 & 4
  \\
   9 & 60 & 1 & . & . & 2 & . & 2 & 2 & 1 & . & 2 & 2 & 1 & 2 & 1 & 1 & 2
   \\  \midrule
  10 & 120 & . & . & . & . & . & . & . & 1 & . & . & . & 1 & . & 1 & 1 & 0
  \\
  10 & 120 & . & 1 & 1 & . & 1 & . & . & . & 1 & . & . & . & . & . & . & 1 
  \\
  10 & 120 & . & 1 & 1 & . & 1 & . & . & . & 1 & . & . & . & . & . & . & 2
  \\
  \bottomrule
  \end{array}
  \]
  \medskip
  \caption{Canonical forms of $2 \times 2 \times 2 \times 2$ symmetric tensors over $\mathbb{F}_5$ (cont.)}
  \label{table2222symmetricmod5-2}
  \end{table}

\subsection{$\mathbb{F}_7$}
There are 33,232,930,569,601 tensors over $\mathbb{F}_7$, 16,807 of which are symmetric. Again, all symmetric tensors have a symmetric decomposition. The maximum symmetric rank in this case is 6, and there are 42 orbits. The number of symmetric tensors in each symmetric rank and the approximate percentages are listed below.
\[
\begin{array}{lrrrrrrr}
\text{rank} & 0 & 1 & 2 & 3 & 4 & 5 & 6 \\
\text{number} & 1 & 24 & 276 & 1932 & 7119 & 6615 & 840 \\
\text{$\approx$ $\%$} & 0.01\% & 0.14\% & 1.64\% & 11.50\% & 42.38\% & 39.36\% & 5.00\%
\end{array}
\]

The symmetric ranks, orders of each orbit, and the minimal representatives of each orbit are given in Table \ref{table2222symmetricmod7}. 

  \begin{table}
  \[
  \begin{array}{ccccccccccccccccccc}
  \text{sym rank} & \text{orbit size} &\multicolumn{16}{c}{\text{canonical form (flattened)}} 
    \\
  \toprule
   0 & 1    & . & . & . & . & . & . & . & . & . & . & . & . & . & . & . & . \\      
  \midrule
   1 &   24    & . & . & . & . & . & . & . & . & . & . & . & . & . & . & . & 1 \\      
  \midrule
   2 &   24    & . & . & . & . & . & . & . & . & . & . & . & . & . & . & . & 3 \\      
   2 &  252    & 1 & . & . & . & . & . & . & . & . & . & . & . & . & . & . & 1 \\      
  \midrule
  3 &   84    & . & . & . & 3 & . & 3 & 3 & . & . & 3 & 3 & . & 3 & . & . & . \\      
  3 &  336    & . & 1 & 1 & . & 1 & . & . & . & 1 & . & . & . & . & . & . & 6 \\      
  3 &  504    & . & 1 & 1 & . & 1 & . & . & 1 & 1 & . & . & 1 & . & 1 & 1 & . \\      
  3 &  504    & . & 1 & 1 & . & 1 & . & . & 3 & 1 & . & . & 3 & . & 3 & 3 & 3 \\      
  3 &  504    & 1 & . & . & 1 & . & 1 & 1 & . & . & 1 & 1 & . & 1 & . & . & 3 \\      
  \midrule
  4 &  336    & . & . & . & . & . & . & . & 1 & . & . & . & 1 & . & 1 & 1 & . \\      
  4 &  504    & . & . & . & 1 & . & 1 & 1 & . & . & 1 & 1 & . & 1 & . & . & 1 \\      
  4 &  504    & . & . & . & 3 & . & 3 & 3 & . & . & 3 & 3 & . & 3 & . & . & 1 \\      
  4 &  336    & . & 1 & 1 & . & 1 & . & . & . & 1 & . & . & . & . & . & . & 1 \\      
  4 &   84    & . & 1 & 1 & . & 1 & . & . & . & 1 & . & . & . & . & . & . & 3 \\      
  4 &  504    & . & 1 & 1 & . & 1 & . & . & 1 & 1 & . & . & 1 & . & 1 & 1 & 2 \\      
  4 & 1008    & . & 1 & 1 & . & 1 & . & . & 1 & 1 & . & . & 1 & . & 1 & 1 & 4 \\      
  4 &  252    & . & 1 & 1 & . & 1 & . & . & 3 & 1 & . & . & 3 & . & 3 & 3 & . \\      
  4 & 1008    & . & 1 & 1 & . & 1 & . & . & 3 & 1 & . & . & 3 & . & 3 & 3 & 1 \\      
  4 &  504    & . & 1 & 1 & . & 1 & . & . & 3 & 1 & . & . & 3 & . & 3 & 3 & 2 \\      
  4 &  504    & . & 1 & 1 & . & 1 & . & . & 3 & 1 & . & . & 3 & . & 3 & 3 & 4 \\      
  4 &  504    & . & 1 & 1 & . & 1 & . & . & 3 & 1 & . & . & 3 & . & 3 & 3 & 5 \\      
  4 &  252    & 1 & . & . & 1 & . & 1 & 1 & . & . & 1 & 1 & . & 1 & . & . & 4 \\      
  4 &  504    & 1 & . & . & 1 & . & 1 & 1 & . & . & 1 & 1 & . & 1 & . & . & 6 \\      
  4 &  252    & 1 & . & . & 1 & . & 1 & 1 & 2 & . & 1 & 1 & 2 & 1 & 2 & 2 & 5 \\      
  4 &   63    & 3 & . & . & 1 & . & 1 & 1 & . & . & 1 & 1 & . & 1 & . & . & 3 \\      
  \midrule
  5 &  504    & . & . & . & 1 & . & 1 & 1 & . & . & 1 & 1 & . & 1 & . & . & 3 \\      
  5 &  504    & . & . & . & 3 & . & 3 & 3 & . & . & 3 & 3 & . & 3 & . & . & 3 \\      
  5 &  336    & . & 1 & 1 & . & 1 & . & . & . & 1 & . & . & . & . & . & . & 2 \\      
  5 &   84    & . & 1 & 1 & . & 1 & . & . & . & 1 & . & . & . & . & . & . & 4 \\      
  5 &  336    & . & 1 & 1 & . & 1 & . & . & . & 1 & . & . & . & . & . & . & 5 \\      
  5 & 1008    & . & 1 & 1 & . & 1 & . & . & 1 & 1 & . & . & 1 & . & 1 & 1 & 3 \\      
  5 &  504    & . & 1 & 1 & . & 1 & . & . & 1 & 1 & . & . & 1 & . & 1 & 1 & 5 \\      
  5 & 1008    & . & 1 & 1 & . & 1 & . & . & 3 & 1 & . & . & 3 & . & 3 & 3 & 6 \\      
  5 &  504    & 1 & . & . & . & . & . & . & 2 & . & . & . & 2 & . & 2 & 2 & 1 \\      
  5 &  504    & 1 & . & . & . & . & . & . & 3 & . & . & . & 3 & . & 3 & 3 & 3 \\      
  5 &  504    & 1 & . & . & 1 & . & 1 & 1 & 1 & . & 1 & 1 & 1 & 1 & 1 & 1 & 2 \\      
  5 &  252    & 1 & . & . & 3 & . & 3 & 3 & . & . & 3 & 3 & . & 3 & . & . & 1 \\      
  5 &   63    & 1 & . & . & 3 & . & 3 & 3 & . & . & 3 & 3 & . & 3 & . & . & 4 \\      
  5 &  504    & 1 & . & . & 3 & . & 3 & 3 & . & . & 3 & 3 & . & 3 & . & . & 5 \\      
 \midrule
  6 &   84    & . & . & . & 1 & . & 1 & 1 & . & . & 1 & 1 & . & 1 & . & . & . \\      
  6 &  504    & 1 & . & . & 1 & . & 1 & 1 & 1 & . & 1 & 1 & 1 & 1 & 1 & 1 & 6 \\      
  6 &  252    & 1 & . & . & 3 & . & 3 & 3 & 3 & . & 3 & 3 & 3 & 3 & 3 & 3 & 5 \\      
  \\
  \bottomrule
  \end{array}
  \]
  \medskip
  \caption{Canonical forms of $2 \times 2 \times 2 \times 2$ symmetric tensors over $\mathbb{F}_7$}
  \label{table2222symmetricmod7}
  \end{table}

We mention now that for larger primes our computer did not have sufficient memory to complete the computations in a reasonable amount of time. 

\section{Conclusion}

In this paper we began by summarizing known results about the canonical forms of $2 \times 2 \times 2$ and $2 \times 2 \times 2 \times 2$ symmetric tensors over $\mathbb{R}$ and $\mathbb{C}$. Our original contribution was considering the same problem over specific prime fields $\mathbb{F}_p$. We obtained a classification of canonical forms of symmetric tensors, and verified that the maximum symmetric rank over certain prime fields is greater than that over $\mathbb{R}$ and $\mathbb{C}$. Based on the prime fields we considered, the maximum symmetric rank of $2 \times 2 \times 2$ symmetric tensors is at least 5. The maximum symmetric rank of $2 \times 2 \times 2 \times 2$ symmetric tensors is at least 10. Due to memory limitations, we could not consider larger prime fields. 

Our algorithms are modifications of those used in \cite{BremnerStavrou}, where we computed canonical forms of $2 \times 2 \times 2$ and $2 \times 2 \times 2 \times 2$ tensors over $\mathbb{F}_2$ and $\mathbb{F}_3$.

\section*{Acknowledgements}
This work is part of my Master's thesis at the University of Saskatchewan. I wish to thank my supervisor, Professor Murray R. Bremner, for the time and valuable input he put into my thesis.


\end{document}